\documentclass[a4paper,10pt]{article}
\usepackage{amssymb}
\usepackage{amsmath,amsfonts,amsthm,amssymb}
\usepackage[dvips]{graphics}
\usepackage{epsfig}
\usepackage{indentfirst}
\usepackage{color}
\allowdisplaybreaks
\newtheoremstyle{theorem}%name
  {10pt}          % space above
  {10pt}  % space below
  {\sl}  % bofy font
 {}% ident - empty=no indent,  \parindent= paragraph indent
  {\bf}  % thm head font
  {. }    % punctuation after thm head
  { }    % space after thm head: `` ``=normal \newline=linebreak
  {}     % thm head specification
\theoremstyle{theorem}

\newtheorem{theorem}{Theorem}[section]

\newtheorem{definition}{Definition}[section]
 \newtheorem{lemma}{Lemma}[section]
 
 \newtheorem{remark}{Remark}[section]
  
\numberwithin{equation}{section}

\newtheoremstyle{defi}%name
  {10pt}          % space above
  {10pt}  % space below
  {\rm}  % bofy font
  {}  % ident - empty=no indent,  \parindent= paragraph indent
  {\bf}  % thm head font
  {. }    % punctuation after thm head
  { }    % space after thm head: `` ``=normal \newline=linebreak
  {}     % thm head specification
\theoremstyle{defi}

\begin{document}
\baselineskip = 13.5pt

\title{\bf On weak-strong uniqueness and singular limit for the compressible Primitive Equations }

\author{ Hongjun Gao$^{1}$  \footnote{Email:gaohj@njnu.edu.cn}\ \ \
\v{S}\'{a}rka Ne\v{c}asov\'{a}$^2$
\footnote{Email: matus@math.cas.cz} \ \ \  Tong Tang$^{3,2}$ \footnote{Email: tt0507010156@126.com}\\
{\small 1.Institute of Mathematics, School of Mathematical Sciences,}\\
{\small Nanjing Normal University, Nanjing 210023, P.R. China}\\
{\small  2. Institute of Mathematics of the Academy of Sciences of the Czech Republic,} \\
{\small \v Zitn\' a 25, 11567, Praha 1, Czech Republic}\\
{\small 3. Department of Mathematics, College of Sciences,}\\
{\small Hohai University, Nanjing 210098, P.R. China}\\
\date{}}

\maketitle
\begin{abstract}
The paper addresses the weak-strong uniqueness property and singular limit for the compressible Primitive Equations (PE). We show that a weak solution coincides with the strong solution emanating from the same initial data. On the other hand, we prove compressible PE will approach to the incompressible inviscid PE equations in the regime of low Mach number and large Reynolds number in the case of well-prepared initial data. To the best of the authors' knowledge, this is the first work to bridge the link between the compressible PE with incompressible inviscid PE.
\vspace{0.5cm}


\medskip

\end{abstract}

\maketitle
\section{Introduction}\setcounter{equation}{0}
The earth is surrounded and occupied by atmosphere and ocean, which play an important role in human's life. From the mathematical point of view and numerical perspective, it is very complicated to use the full hydrodynamical and thermodynamical equations to describe the motion and fascinating phenomena of atmosphere and ocean. In order to simplify model, scientists introduce the Primitive Equations (PE) model in meteorology and geophysical fluid dynamics, which helps us to predict the long-term weather and detect the global climate changes. In this paper, we study the following Compressible Primitive Equations (CPE):

\begin{eqnarray}
\left\{
\begin{array}{llll}  \partial_{t}\rho+\text{div}_x(\rho \mathbf{u})+\partial_z(\rho w)=0, \\
\partial_t(\rho \mathbf{u})+\textrm{div}_x(\rho\mathbf{u}\otimes\mathbf{u})+\partial_z(\rho\mathbf uw)+\nabla_x p(\rho)=\mu\Delta_x\mathbf u+\lambda\partial^2_{zz}\mathbf u,\\
\partial_zp(\rho)=0,
\end{array}\right.\label{a}
\end{eqnarray}
in $(0,T)\times\Omega$. Here $\Omega=\{(x,z)|x\in\mathbb{T}^2,0<z<1\}$, $x$ denotes the horizontal direction and $z$ denotes the vertical direction.  $\rho=\rho(t,x)$, $\mathbf{u}(t,x,z)\in\mathbb{R}^2$ and $w(t,x,z)\in\mathbb{R}$
represent the density, the horizonal velocity and vertical velocity respectively.

From the hydrostatic balance equation $(1.1)_3$, it follows that {\bf the density $\rho$ is independent of $z$}.
$\mu>0$, $\lambda\geq0$ are the constant viscosity coefficients. The system is supplemented by the boundary conditions
\begin{eqnarray}
w|_{z=0}=w|_{z=1}=0,\hspace{4pt}\partial_z\mathbf u|_{z=0}=\partial_z\mathbf u|_{z=1}=0,
\end{eqnarray}
and initial data
\begin{eqnarray}
\rho\mathbf u|_{t=0}=\mathbf m_0(x,z),\hspace{3pt}\rho|_{t=0}=\rho_0(x).
\end{eqnarray}

The pressure $p(\rho)$ satisfies the barotropic pressure law where the pressure and the density are related by the following formula:
\begin{eqnarray}
p(\rho)=\rho^\gamma\hspace{5pt}(\gamma>1).
\end{eqnarray}
The PE model is widely used in meteorology and geophysical fluid dynamics, due to its accurate theoretical analysis and practical numerical computing. Concerning geophysical fluid dynamics we can refer to work by Chemin, Desjardins, Gallagher and Grenier \cite{ch} or Feireisl, Gallagher, Novotn\'{y} \cite{e}. There is a great number of results about PE, such as \cite{bg,b1,c2,c4,c5,l3,l4,s,t,ws}. We just mention some of results. Guill\'{e}n-Gonz\'{a}lez, Masmoudi and Rodr\'{\i}guez-Bellido \cite{gu} proved the local existence of strong solutions. The celebrated breakthrough result was made by Cao and Titi \cite{c1}. They were first who proved the global well-posedness of PE. After that a lot of scientists  were focused on the dynamics and regularity of PE e.g. \cite{g1,g2,ju,kukavica}. Recently in \cite{c2,c4,c5}, the authors considered the strong solution for PE with vertical eddy diffusivity and only horizontal dissipation. About random perturbations of PE, the local and global strong
solution of PE can be referred to \cite{d1, d2, gao}, large deviation principles, see \cite{dong} and diffusion limit, see \cite{g3}. On the other hand, regarding to inviscid PE (hydrostatic incompressible Euler equations), the existence and uniqueness is an outstanding open problem. Only a few  results are available. Under the convex horizontal velocity assumptions, Brenier \cite{b} proved the existence of smooth solutions in two-dimensions. Then, Masmoudi and Wong \cite{m} utilized the weighted $H^s$ a priori estimates and obtained the existence, uniqueness and weak-strong uniqueness. Removing the convex horizontal velocity assumptions, they extended Brenier's result. By virtue of Cauchy-Kowalewshi theorem, the authors \cite{k} constructed a locally, unique and real-analytic solution. Notably, Brenier \cite{by} suggested that the existence problem may be ill-posed in Sobolev spaces. Further Cao et al. \cite{c3} established the blow up for certain class of smooth solutions in finite time.

In order to show the atmosphere and ocean have compressible property, Ersoy et al. \cite{er1} consider that the vertical scale of atmosphere is significantly smaller than the horizontal scales and  they derive the CPE from the compressible Navier-Stokes equations. To be precise, the CPE system is obtained by replacing the vertical velocity momentum equation with hydrostatic balance equation. Compared with compressible Navier-Stokes equations, the regularity of vertical velocity is less regular than horizontal velocity in the CPE system. In the absence of sufficient information about the vertical velocity, it inevitably leads to difficulty for obtaining the existence of solutions.  \emph{ Lions, Teman and Wang \cite{l1,l2} were first to  study the CPE and received fundamental results in this field.} Under a smart $P-coordinates$, they reformulated the system into the classical PE with the incompressible condition. Later on, Gatapov and Kazhikhov \cite{g}, Ersoy and Ngom \cite{er2} proved the global existence of weak solutions in 2D case. Liu and Titi \cite{liu1} used the classical methods to proved the local existence of strong solutions in 3D case. Ersoy et al. \cite{er1}, Tang and Gao \cite{tang} showed the stability of weak solutions with the viscosity coefficients depending on the density. The stability means that a subsequence of weak solutions will converge to another weak solutions if it satisfies some uniform bounds. Recently, based on the work \cite{b1,b2,b3,li,v},  Liu and Titi \cite{liu2} and independently  Wang et al. \cite{w} used the B-D entropy to prove the global existence of weak solutions in the case where the viscosity coefficients are depending on the density.

Our paper is divided into two parts. The first part concerns the weak-strong uniqueness of CPE. Recently, Liu and Titi \cite{liu3} studied the zero Mach number limit of CPE, proving it converges to incompressible PE, which is a breakthrough result to bridge the link between CPE and PE system. In the second part, inspired by \cite{liu3}, we investigate the singular limit of CPE, showing it converges to incompressible inviscid PE system. \emph{This is the first attempt to use the relative entropy method to study asymptotic limit for CPE.} Let us mention that the corner-stone analysis of our results is based on the relative energy inequality which was invented by Dafermos, see \cite{D}. It was introduced by  Germain \cite{ge} and generalized by Feireisl \cite{e2} for compressible fluid model. Feireisl and his  co-authors \cite{e3,e4} used the versatile tool to solve various problems. However, compared with the previous classical results, there is significant difference in the process of using relative energy inequality to CPE model due to the absence of the information on the vertical velocity. Therefore, it is not straightforward to apply the method from Navier-Stokes to CPE. We utilize the special structure of CPE to find the deeper relationship and reveal the important feature of CPE.

The paper is organized as follows. In Section 2, we introduce the dissipative weak solutions, relative energy and state our first theorem. In Section 3, we prove the weak-strong uniqueness. We recall the target system, state the singular limit theorem and derive the necessary uniform bounds in Section 4. Section 5 is devoted to proof of the convergence in the case of well-prepared initial data.

\vskip 0.5cm
%%%%%%%%%%%%%%%%%%%%%%%%%%%%%%%%%%%%%%%%%%%%%%%%%%%%%%%%%%%%%%%%
%%%%%%%%%%%%%%%%%%%%%%%%%%%%%%%%%%%%%%%%%%%%%%%%%%%%%%%%%%%%%%%%
\noindent {\bf Part I: Weak-Strong uniqueness}
\vskip 0.2cm
In this part,  we focus on the weak-strong uniqueness of the CPE system.
\section{Preliminaries and main result}

First of all, we should point out that a proper notion of weak solution to CPE has not been well understood. Recently, Bresch and Jabin \cite{br} consider different compactness method from Lions or Feireisl which can be applied to anisotropical stress tensor similarly. They  obtain the global existence of weak solutions if $|\mu-\lambda|$ are not too large. Let us state one of the possible definitions here.

\subsection{Dissipative weak solutions}

\begin{definition}\label{def1}
We say that $[\rho,\mathbf u,w]$ is a dissipative weak solution to the system of \eqref{a}, supplemented with initial data (1.3) and pressure follows (1.4) if $\rho=\rho(x,t)$ and
\begin{align}
\mathbf u\in L^2(0,T;H^1(\Omega)),\hspace{3pt} \rho|\mathbf u|^2\in L^\infty(0,T; L^1(\Omega)),
\hspace{3pt}\rho\in L^\infty(0,T;L^\gamma(\Omega)\cap L^1(\Omega)).
\end{align}

\noindent
$\bullet$ the continuity equation
\begin{align}
[\int_\Omega\rho\varphi dxdz]^{t=\tau}_{t=0}=\int^\tau_0\int_{\Omega}\rho\partial_t\varphi+\rho\mathbf{u}\nabla_x\varphi+\rho w\partial_z\varphi dxdzdt,
\end{align}
holds for all $\varphi\in C^\infty_c([0,T)\times\Omega)$;

\noindent
$\bullet$
the momentum equation
\begin{align}
[\int_\Omega\rho\mathbf u\varphi dxdz]^{t=\tau}_{t=0}&=\int^\tau_0\int_{\Omega}\rho\mathbf{u}\partial_t\varphi+ \rho\mathbf{u}\otimes\mathbf{u}:\nabla_x\varphi+\rho\mathbf uw\partial_z\varphi+ p(\rho)\text{div}\varphi dxdzdt\nonumber\\
&\hspace{8pt}-\int^\tau_0\int_{\Omega}[\mu\nabla_x\mathbf u:\nabla_x\varphi+\lambda\partial_z\mathbf u\partial_z\varphi]dxdzdt,
\end{align}
holds for all $\varphi\in C^\infty_c([0,T)\times\Omega)$,

\noindent
$\bullet$
the energy inequality
\begin{align}
[\int_{\Omega}\frac{1}{2}\rho|\mathbf{u}|^2+P(\rho)-P'(\overline{\rho})(\rho-\overline{\rho})-P(\overline{\rho})dxdz]|^{t=\tau}_{t=0}
+\int^\tau_0\int_\Omega(\mu|\nabla_x\mathbf u|^2+\lambda|\partial_z\mathbf u|^2)dxdzdt\leq 0,
\end{align}
holds for a.a $\tau\in(0,T)$, a arbitrary constant $\overline{\rho}$, where $P(\rho)=\rho\int^\rho_1\frac{p(z)}{z^2}dz$.

Moreover, as there is no information about $w$, so we need the following equation:
\begin{align}
\rho w(x,z,t)=-\rm{div}_x(\rho\widetilde{\mathbf u})+z\rm{div}_x(\rho\overline{\mathbf u}), \hspace{4pt}
\text{in the sense of} \hspace{4pt}H^{-1}(\Omega),
\label{b1}
\end{align}
where
\begin{align*}
\widetilde{\mathbf u}(x,z,t)=\int^z_0u(x,s,t)ds,\hspace{5pt}\overline{u}(x,t)=\int^1_0u(x,z,t)dz.
\end{align*}

\end{definition}
We should emphasize that \eqref{b1} is the key step to obtain the existence of weak solution in \cite{liu2,w}, which is inspired by incompressible case.

\subsection{Relative entropy inequality}

Motivated by \cite{e2,e3}, for any finite weak solution $(\rho,\mathbf u,w)$ to the CPE system, we introduce the relative energy functional
\begin{align}
\mathcal{E}(\rho,\mathbf{u}|r, \mathbf{U})&=\int_{\Omega}[\frac{1}{2}\rho|\mathbf u-\mathbf U|^2+P(\rho)-P'(r)(\rho-r)-P(r)]dxdz\nonumber\\
&=\int_\Omega(\frac{1}{2}\rho^2|\mathbf u|+P(\rho))dxdz-\int_\Omega\rho\mathbf u\cdot\mathbf Udxdz
+\int_\Omega\rho[\frac{|\mathbf U|^2}{2}-P'(r)]dxdz+\int_\Omega p(r)dxdz\nonumber\\
&=\sum^4_{i=1}I_i,\label{a1}
\end{align}
where $r>0$, $\mathbf U$ are smooth ``test'' functions, $r$, $\mathbf U$ compactly supported in $\Omega$.

\begin{lemma}\label{relativeentropy}
Let $(\rho,\mathbf{u}, w)$ be a dissipative weak solution introduced in Definition  \ref{def1}. Then $(\rho,\mathbf{u}, w)$
satisfy the relative entropy inequality

\begin{align}
\mathcal{E}&(\rho,\mathbf{u}|r,\mathbf U)|^{t=\tau}_{t=0}+\int^\tau_0\int_\Omega\big{(}\mu\nabla_x\mathbf u\cdot(\nabla_x\mathbf u-\nabla_x\mathbf U)+\lambda\partial_z\mathbf u(\partial_z\mathbf u-\partial_z\mathbf U)\big{)}dxdzdt\nonumber\\
&\leq\int^\tau_0\int_{\Omega}\rho(\mathbf U-\mathbf u)\partial_t\mathbf U+\rho\mathbf u(\mathbf U-\mathbf u)\cdot\nabla_x\mathbf U
+\rho w(\mathbf U-\mathbf u)\cdot\partial_z\mathbf U-p(\rho)\text{div}_x\mathbf Udxdzdt\nonumber\\
&\hspace{15pt}-\int^\tau_0\int_{\Omega}P''(r)(\rho\partial_tr+\rho\mathbf u\nabla_xr)dxdzdt
+\int^\tau_0\int_{\Omega}\partial_tp(r)dxdzdt.
\end{align}
\end{lemma}
{\bf Proof:}
From the weak formulation and energy inequality (2.2)-(2.4) we deduce
\begin{align}
&I_1|^{t=\tau}_{t=0}+\int^\tau_0\int_\Omega(\mu|\nabla_x\mathbf u|^2+\lambda|\partial_z\mathbf u|^2)dxdzdt\leq0,\\
&I_2|^{t=\tau}_{t=0}=-\int^\tau_0\int_\Omega\rho\mathbf u\partial_t\mathbf U+\rho\mathbf u\otimes\mathbf u:\nabla_x\mathbf U+
\rho\mathbf uw\partial_z\mathbf U+p(\rho)\text{div}_x\mathbf Udxdzdt\nonumber\\
&\hspace{40pt}+\int^\tau_0\int_\Omega\mu\nabla_x\mathbf u:\nabla_x\mathbf U+\lambda\partial_z\mathbf u\partial_z\mathbf Udxdzdt,\\
&I_3|^{t=\tau}_{t=0}=\int^\tau_0\int_\Omega\rho\partial_t\frac{|\mathbf U|^2}{2}+\rho\mathbf u\cdot\nabla_x\frac{|\mathbf U|^2}{2}+\rho w\partial_z\frac{|\mathbf U|^2}{2}dxdzdt\nonumber\\
&\hspace{40pt}-\int^\tau_0\int_\Omega\rho\partial_tP'(r)+\rho\mathbf u\cdot\nabla_xP'(r)+\rho w\partial_zP'(r)dxdzdt\nonumber\\
&\hspace{20pt}=\int^\tau_0\int_\Omega\rho\mathbf U\partial_t\mathbf U+\rho\mathbf u\mathbf U\cdot\nabla_x\mathbf U+\rho w\mathbf U\partial_z\mathbf Udxdzdt\nonumber\\
&\hspace{30pt}-\int^\tau_0\int_\Omega\rho P''(r)\partial_tr+P''(r)\rho\mathbf u\cdot\nabla_xr dxdzdt,\\
&I_4|^{t=\tau}_{t=0}=[\int_\Omega p(\rho)dxdz]|^{t=\tau}_{t=0}.
\end{align}

Summing (2.6)-(2.10) together, we obtain
\begin{align}
\mathcal{E}&(\rho,\mathbf{u}|r,\mathbf U)|^{t=\tau}_{t=0}+\int^\tau_0\int_\Omega\big{(}\mu\nabla_x\mathbf u\cdot(\nabla_x\mathbf u-\nabla_x\mathbf U)+\lambda\partial_z\mathbf u(\partial_z\mathbf u-\partial_z\mathbf U)\big{)}dxdzdt\nonumber\\
&\leq\int^\tau_0\int_{\Omega}\rho(\mathbf U-\mathbf u)\partial_t\mathbf U+\rho\mathbf u(\mathbf U-\mathbf u)\cdot\nabla_x\mathbf U
+\rho w(\mathbf U-\mathbf u)\cdot\partial_z\mathbf U-p(\rho)\text{div}_x\mathbf Udxdzdt\nonumber\\
&\hspace{15pt}-\int^\tau_0\int_{\Omega}P''(r)(\rho\partial_tr+\rho\mathbf u\nabla_xr)dxdzdt
+\int^\tau_0\int_{\Omega}\partial_tp(r)dxdzdt.
\end{align}
%%%%%%%%%%%%%%%%%%%%%%%%%%%%%%%%%%%%%%%%%%%%%%%%%%%%%%%%%%%%%%%%%%%
%%%%%%%%%%%%%%%%%%%%%%%%%%%%%%%%%%%%%%%%%%%%%%%%%%%%%%%%%%%%%%%%
%%%%%%%%%%%%%%%%%%%%%%%%%%%%%%%%%%%%%%%%%%%%%%%%%%%%%%%%%%%%%%%%
\subsection{Main result}
We say that $(r,\mathbf U,W)$ is a strong solution to the CPE system $(1.1)-(1.4)$ in $(0,T)\times\Omega$, if
\begin{align*}
&r^\frac{1}{2}\in L^\infty(0,T;H^2(\Omega)),\hspace{3pt}\partial_tr^\frac{1}{2}\in L^\infty(0,T;H^1(\Omega)),\hspace{3pt}r>0\hspace{3pt}\text{for all}\hspace{3pt}(t,x),\\
&\mathbf U\in L^\infty(0,T;H^3(\Omega))\cap L^2(0,T;H^4(\Omega)),\hspace{3pt} \partial_t\mathbf U\in L^2(0,T; H^2(\Omega)),
\end{align*}
with initial data $r^\frac{1}{2}_0\in H^2(\Omega)$, $r_0>0$ and $\mathbf U_0\in H^3(\Omega)$.

Now, we are ready to state our first result.
\begin{theorem}
Let $\gamma>6$, $(\rho,\mathbf u,w)$ be a dissipative weak solution to the CPE system $(1.1)-(1.4)$ in $(0,T)\times\Omega$. Let $(r,\mathbf U, W)$ be a strong solution to the same problem and emanating from the same initial data. Then,
\begin{align*}
\rho=r,\hspace{5pt}\mathbf u=\mathbf U,\hspace{4pt}\text{in}\hspace{3pt}(0,T)\times\Omega.
\end{align*}

\end{theorem}

\begin{remark}
Liu and Titi \cite{liu1} obtained the local existence of strong solutions to CPE. It is important to point out that our result holds under more regularity than the strong solutions obtained in \cite{liu1}.
\end{remark}

Section 3 is devoted to the proof of the above theorem.

%%%%%%%%%%%%%%%%%%%%%%%%%%%%%%%%%%%%%%%
\section{Weak-strong uniqueness}
The proof of Theorem 2.1 depends on the relative energy inequality by considering the strong solution $[r,\mathbf U, W]$ as test function in the relative energy inequality \eqref{a1}.

\subsection{Step 1}
We write
\begin{align*}
\int_\Omega\rho\mathbf u(\mathbf U-\mathbf u)\cdot\nabla_x\mathbf Udxdz=
\int_\Omega\rho(\mathbf u-\mathbf U)(\mathbf U-\mathbf u)\cdot\nabla_x\mathbf Udxdz
+\int_\Omega\rho\mathbf U(\mathbf U-\mathbf u)\cdot\nabla_x\mathbf Udxdz.
\end{align*}

As $[r,\mathbf U, W]$ is a strong solution, it is easy to obtain that
\begin{align}
\int_\Omega\rho(\mathbf u-\mathbf U)(\mathbf U-\mathbf u)\cdot\nabla_x\mathbf Udxdz
\leq C\mathcal{E}(\rho,\mathbf u|r,\mathbf U).
\end{align}

Moreover, the momentum equation reads as
\begin{align*}
(r\mathbf U)_t+\text{div}(r\mathbf U\otimes\mathbf U)+\partial_z(r\mathbf UW)+\nabla_xp(r)=\mu\Delta_x\mathbf U+\lambda\partial_{zz}\mathbf U,
\end{align*}
implying that
\begin{align*}
\mathbf U_t+\mathbf U\cdot\nabla_x\mathbf U+W\partial_z\mathbf U=-\frac{1}{r}\nabla_xp(r)+\frac{\mu}{r}\Delta_x\mathbf U
+\frac{\lambda}{r}\partial_{zz}\mathbf U.
\end{align*}

So we rewrite
\begin{align*}
\int_\Omega\rho(\mathbf U-\mathbf u)\cdot\partial_t\mathbf U+
\rho\mathbf U(\mathbf U-\mathbf u)\cdot\nabla_x\mathbf U
+\rho W(\mathbf U-\mathbf u)\cdot\partial_z\mathbf U+\rho(w-W)(\mathbf U-\mathbf u)\cdot\partial_z\mathbf Udxdz\\
=\int_\Omega\frac{\rho}{r}(\mathbf U-\mathbf u)(-\nabla_xp(r)+\mu\Delta_x\mathbf U+\lambda\partial_{zz}\mathbf U)dxdz
+\int_\Omega\rho(w-W)(\mathbf U-\mathbf u)\cdot\partial_z\mathbf Udxdz.
\end{align*}

Thus, we obtain that
\begin{align*}
\mathcal{E}&(\rho,\mathbf{u}|r,\mathbf U)|^{t=\tau}_{t=0}+\int^\tau_0\int_\Omega\big{(}\mu\nabla_x\mathbf u\cdot(\nabla_x\mathbf u-\nabla_x\mathbf U)+\lambda\partial_z\mathbf u(\partial_z\mathbf u-\partial_z\mathbf U)\big{)}dxdzdt\nonumber\\
&\leq C\int^\tau_0\mathcal{E}(\rho,\mathbf{u}|r,\mathbf U)dt
-\int^\tau_0\int_{\Omega}P''(r)(\rho\partial_tr+\rho\mathbf u\nabla_xr)dxdzdt\nonumber\\
&\hspace{8pt}+\int^\tau_0\int_\Omega\frac{\rho}{r}(\mathbf U-\mathbf u)(\mu\Delta_x\mathbf U+\lambda\partial_{zz}\mathbf U)dxdz
-\int^\tau_0\int_\Omega\frac{\rho}{r}(\mathbf U-\mathbf u)\nabla_xp(r)dxdz\nonumber\\
&\hspace{8pt}+\int^\tau_0\int_\Omega\rho(w-W)(\mathbf U-\mathbf u)\cdot\partial_z\mathbf Udxdzdt+\int^\tau_0\int_\Omega\partial_tp(r)dxdzdt
-\int^\tau_0\int_\Omega p(\rho)\text{div}_x\mathbf Udxdzdt.\nonumber\\
\end{align*}

Before estimating, we should recall the following useful inequality from \cite{e2}:
\begin{equation} \label{pres}
P(\rho)-P'(r)(\rho-r)-P(r)\geq\left\{
\begin{array}{llll} C|\rho-r|^2,\hspace{5pt}\text{when} \hspace{3pt} \frac{r}{2}<\rho<r, \nonumber\\
C(1+\rho^\gamma),\hspace{5pt}\text{otherwise}.
\end{array}\right.
\end{equation}

Moreover, from \cite{e2}, we learn that
\begin{align}
&\mathcal{E}(\rho,\mathbf{u}|r,\mathbf U)(t)\in L^\infty(0,T),\hspace{3pt}
\int_\Omega\chi_{\rho\geq r}\rho^{\gamma}dxdz\leq C\mathcal{E}(\rho,\mathbf{u}|r,\mathbf U)(t),\nonumber\\
&\int_\Omega\chi_{\rho\leq \frac{r}{2}}1dxdz\leq C\mathcal{E}(\rho,\mathbf{u}|r,\mathbf U)(t),\hspace{3pt}
\int_\Omega\chi_{\frac{r}{2}<\rho<r}(\rho-r)^2dxdz\leq C\mathcal{E}(\rho,\mathbf{u}|r,\mathbf U)(t).\label{a3}
\end{align}

The main difficulty is to estimate the complicated nonlinear term $\int_\Omega\rho(w-W)(\mathbf U-\mathbf u)\cdot\partial_z\mathbf Udxdz$, we rewrite it as
\begin{align}
\int_\Omega\rho(w-W)(\mathbf U-\mathbf u)\cdot\partial_z\mathbf Udxdz
=\int_\Omega\rho w(\mathbf U-\mathbf u)\cdot\partial_z\mathbf Udxdz-\int_\Omega\rho W(\mathbf U-\mathbf u)\cdot\partial_z\mathbf Udxdz.\label{b}
\end{align}

According to \cite{e2,kr}, we divide the second term on the right side of (3.3) into three parts
\begin{align}
\int_\Omega&\rho W(\mathbf U-\mathbf u)\cdot\partial_z\mathbf Udxdz\nonumber\nonumber\\
&=\int_\Omega\chi_{\rho\leq \frac{r}{2}}\rho W(\mathbf U-\mathbf u)\cdot\partial_z\mathbf Udxdz
+\int_\Omega\chi_{\frac{r}{2}<\rho<r}\rho W(\mathbf U-\mathbf u)\cdot\partial_z\mathbf Udxdz+\int_\Omega\chi_{\rho\geq r}\rho W(\mathbf U-\mathbf u)\cdot\partial_z\mathbf Udxdz\nonumber\nonumber\\
&\leq \|\chi_{\rho\leq \frac{r}{2}}1\|_{L^2(\Omega)}\|r\|_{L^\infty}\|W\partial_z\mathbf U\|_{L^3}\|\mathbf U-\mathbf u\|_{L^6(\Omega)}
+\int_\Omega\chi_{\rho\geq r}\rho^{\frac{\gamma}{2}}W\partial_z\mathbf U\cdot(\mathbf U-\mathbf u)dxdz\nonumber\\
&\hspace{8pt}+C\|\chi_{\frac{r}{2}<\rho< r}(\rho-r)\|_{L^2(\Omega)}\|W\partial_z\mathbf U\|_{L^3}\|\mathbf U-\mathbf u\|_{L^6(\Omega)}\nonumber\\
&\leq C\int_\Omega\chi_{\rho\leq \frac{r}{2}}1dxdz +C\int_\Omega\chi_{\frac{r}{2}<\rho<r}(\rho-r)^2dxdz
+C\int_\Omega\chi_{\rho\geq r}\rho^\gamma dxdz
+\delta\|\mathbf U-\mathbf u\|^2_{L^6(\Omega)}\nonumber\\
&\leq C\mathcal{E}(\rho,\mathbf u|r,\mathbf U)+\delta\|\nabla_x\mathbf U-\nabla_x\mathbf u\|^2_{L^2(\Omega)}
+\delta\|\partial_z\mathbf U-\partial_z\mathbf u\|^2_{L^2(\Omega)},\label{33c}
\end{align}
where in the last inequality, we have used the following celebrated inequality from Feireisl \cite{e1}:
\begin{lemma}
Let $2\leq p\leq6$, and $\rho\geq0$ such that $0<\int_\Omega\rho dx\leq M$ and $\int_\Omega\rho^\gamma dx\leq E_0$ for some $(\gamma>1)$ then
\begin{align*}
\|f\|_{L^p(\Omega)}\leq C\|\nabla f\|_{L^2(\Omega)}+\|\rho^{\frac{1}{2}}f\|_{L^2(\Omega)},
\end{align*}
where $C$ depends on $M$ and $E_0$.
\end{lemma}

On the other hand, we take \eqref{b1} into the first term on the right hand of \eqref{b} and get
\begin{align}
\int_\Omega\rho w&(\mathbf U-\mathbf u)\cdot\partial_z\mathbf Udxdz\nonumber\\
&=\int_\Omega[-\text{div}_x(\rho\widetilde{\mathbf u})+z\text{div}_x(\rho\overline{\mathbf u})]\partial_z\mathbf U
\cdot(\mathbf U-\mathbf u)dxdz\nonumber\\
&=\int_\Omega(\rho\widetilde{\mathbf u}-z\rho\overline{\mathbf u})\partial_z\nabla_x\mathbf U\cdot(\mathbf U-\mathbf u)dxdz+\int_\Omega(\rho\widetilde{\mathbf u}-z\rho\overline{\mathbf u})\partial_z\mathbf U\cdot(\nabla_x\mathbf U-\nabla_x\mathbf u)dxdz.\label{c}
\end{align}

In the following, we will estimate the terms on the right hand side of \eqref{c}. We choose the most complicated terms as examples to estimate, the remaining terms can be analyzed similarly. Firstly, we deal with $\int_\Omega\rho\widetilde{\mathbf u}\partial_z\nabla_x\mathbf U\cdot(\mathbf U-\mathbf u)dxdz$ in the following,
\begin{align*}
\int_\Omega\rho\widetilde{\mathbf u}\partial_z\nabla_x\mathbf U\cdot(\mathbf U-\mathbf u)dxdz
&=\int_\Omega\rho(\widetilde{\mathbf u}-\widetilde{\mathbf U})\partial_z\nabla_x\mathbf U\cdot(\mathbf U-\mathbf u)dxdz
+\int_\Omega\rho\widetilde{\mathbf U}\partial_z\nabla_x\mathbf U\cdot(\mathbf U-\mathbf u)dxdz\\
&=J_1+J_2.
\end{align*}
where $\widetilde{\mathbf U}=\int^z_0\mathbf U(x,s,t)ds$.

Similar to the above analysis, we divide the term $J_2$ into three parts
\begin{align*}
J_2&=\int_\Omega\rho\widetilde{\mathbf U}\partial_z\nabla_x\mathbf U\cdot(\mathbf U-\mathbf u)dxdz\\
&=\int_\Omega\chi_{\rho\leq \frac{r}{2}}\rho\widetilde{\mathbf U}\partial_z\nabla_x\mathbf U\cdot(\mathbf U-\mathbf u)dxdz
+\int_\Omega\chi_{\frac{r}{2}<\rho<r}\rho\widetilde{\mathbf U}\partial_z\nabla_x\mathbf U\cdot(\mathbf U-\mathbf u)dxdz\\
&\hspace{5pt}+\int_\Omega\chi_{\rho\geq r}\rho\widetilde{\mathbf U}\partial_z\nabla_x\mathbf U\cdot(\mathbf U-\mathbf u)dxdz\\
&\leq \|\chi_{\rho\leq \frac{r}{2}}1\|_{L^2(\Omega)}\|r\|_{L^\infty}\|\widetilde{\mathbf U}\partial_z\nabla_x\mathbf U\|_{L^3}\|\mathbf U-\mathbf u\|_{L^6(\Omega)}
+\|\chi_{\rho\geq r}\rho^{\frac{\gamma}{2}}\|_{L^2(\Omega)}\|\widetilde{\mathbf U}\partial_z\nabla_x\mathbf U\|_{L^3(\Omega)}\|\mathbf U-\mathbf u\|_{L^6(\Omega)}\\
&\hspace{8pt}+C\|\chi_{\frac{r}{2}<\rho<r}(\rho-r)\|_{L^2(\Omega)}\|\widetilde{\mathbf U}\partial_z\nabla_x\mathbf U\|_{L^3(\Omega)}\|\mathbf U-\mathbf u\|_{L^6(\Omega)}\\
&\leq C\mathcal{E}(\rho,\mathbf u|r,\mathbf U)(t)+\delta\|\nabla_x\mathbf U-\nabla_x\mathbf u\|^2_{L^2(\Omega)}
+\delta\|\partial_z\mathbf U-\partial_z\mathbf u\|^2_{L^2(\Omega)}.
\end{align*}

On the other hand, by virtue of Cauchy inequality, we obtain
\begin{align}
J_1&=\int_\Omega\rho(\widetilde{\mathbf u}-\widetilde{\mathbf U})\partial_z\nabla_x\mathbf U\cdot(\mathbf U-\mathbf u)dxdz\nonumber\\
&\leq\|\partial_z\nabla_x\mathbf U\|_{L^\infty}\int_\Omega\rho|\widetilde{\mathbf u}-\widetilde{\mathbf U}|^2dxdz+\int_\Omega\rho|\mathbf u-\mathbf U|^2dxdz\nonumber\\
&\leq C\int_\Omega\rho|\int^z_0(\mathbf u(s)-\mathbf U(s))ds|^2dxdz+\mathcal{E}(\rho,\mathbf u|r,\mathbf U)\nonumber\\
&\leq C\int_\Omega\rho\big{(}\int^1_0|\mathbf u-\mathbf U|^2ds\big{)}dxdz+\mathcal{E}(\rho,\mathbf u|r,\mathbf U)\nonumber\\
&\leq C\int^1_0\int_\Omega \rho|\mathbf u-\mathbf U|^2dxdzds+\mathcal{E}(\rho,\mathbf u|r,\mathbf U)\nonumber\\
&\leq C\int_\Omega \rho|\mathbf u-\mathbf U|^2dxdz+\mathcal{E}(\rho,\mathbf u|r,\mathbf U)\nonumber\\
&\leq C\mathcal{E}(\rho,\mathbf u|r,\mathbf U).\label{3aaa}
\end{align}

Secondly, we will tackle with another complicated nonlinear term $\int_\Omega\rho\widetilde{\mathbf u}\partial_z\mathbf U\cdot(\nabla_x\mathbf U-\nabla_x\mathbf u)dxdz$. It is easy to rewrite it as
\begin{align}
\int_\Omega\rho\widetilde{\mathbf u}&\partial_z\mathbf U\cdot(\nabla_x\mathbf U-\nabla_x\mathbf u)dxdz\nonumber\\
&=\int_\Omega\chi_{\rho< r}\rho\widetilde{\mathbf u}\partial_z\mathbf U\cdot(\nabla_x\mathbf U-\nabla_x\mathbf u)dxdz
+\int_\Omega\chi_{\rho\geq r}\rho\widetilde{\mathbf u}\partial_z\mathbf U\cdot(\nabla_x\mathbf U-\nabla_x\mathbf u)dxdz,\label{2a}
\end{align}
where
\begin{align*}
\int_\Omega\chi_{\rho< r}&\rho\widetilde{\mathbf u}\partial_z\mathbf U\cdot(\nabla_x\mathbf U-\nabla_x\mathbf u)dxdz\\
&=\int_{\Omega}\chi_{\rho<r}\rho(\widetilde{\mathbf u}-\widetilde{\mathbf U})\partial_z\mathbf U\cdot(\nabla_x\mathbf U-\nabla_x\mathbf u)dxdz
+\int_{\Omega}\chi_{\rho<r}\rho\widetilde{\mathbf U}\partial_z\mathbf U\cdot(\nabla_x\mathbf U-\nabla_x\mathbf u)dxdz\\
&=\int_{\Omega}\chi_{\rho<r}\rho(\widetilde{\mathbf u}-\widetilde{\mathbf U})\partial_z\mathbf U\cdot(\nabla_x\mathbf U-\nabla_x\mathbf u)dxdz
+\int_{\Omega}\chi_{\frac{r}{2}<\rho<r}\rho\widetilde{\mathbf U}\partial_z\mathbf U\cdot(\nabla_x\mathbf U-\nabla_x\mathbf u)dxdz\\
&\hspace{20pt}+\int_{\Omega}\chi_{\rho\leq \frac{r}{2}}\rho\widetilde{\mathbf U}\partial_z\mathbf U\cdot(\nabla_x\mathbf U-\nabla_x\mathbf u)dxdz\\
&\leq \|\chi_{\rho<r}\rho^{\frac{1}{2}}\|_{L^\infty(\Omega)}\|\sqrt{\rho}(\widetilde{\mathbf u}-\widetilde{\mathbf U})\|_{L^2(\Omega)}\|\partial_z\mathbf U\|_{L^\infty(\Omega)}\|\nabla_x\mathbf U-\nabla_x\mathbf u\|_{L^2(\Omega)}\\
&\hspace{10pt}+\|\chi_{\frac{r}{2}<\rho< r}\rho\|_{L^2(\Omega)}\|\widetilde{\mathbf U}\partial_z\mathbf U\|_{L^\infty(\Omega)}\|\nabla_x\mathbf U-\nabla_x\mathbf u\|_{L^2(\Omega)}\\
&\hspace{10pt}+\|\chi_{\rho\leq \frac{r}{2}}1\|_{L^2(\Omega)}\|r\|_{L^\infty(\Omega)}
\|\widetilde{\mathbf U}\partial_z\mathbf U\|_{L^\infty(\Omega)}\|\nabla_x\mathbf U-\nabla_x\mathbf u\|_{L^2(\Omega)}\\
&\leq C\mathcal{E}(\rho,\mathbf u|r,\mathbf U)(t)+\delta\|\nabla_x\mathbf U-\nabla_x\mathbf u\|^2_{L^2(\Omega)}.
\end{align*}

%From the G-N inequality
%\begin{align*}
%\|D^ju\|_{L^p(\Omega)}\leq \|D^mu\|^\theta_{L^r}\|u\|^{1-\theta}_{L^q},\hspace{5pt}
%\frac{1}{p}=\frac{j}{N}+(\frac{1}{r}-\frac{m}{N})\theta+\frac{1-\theta}{q}
%\end{align*}

Then we will deal with the second term on the right side of \eqref{2a}:
\begin{align}
\int_\Omega\chi_{\rho\geq r}&\rho\widetilde{\mathbf u}\partial_z\mathbf U\cdot(\nabla_x\mathbf U-\nabla_x\mathbf u)dxdz\nonumber\\
&=\int_\Omega\chi_{\rho\geq r}\rho(\widetilde{\mathbf u}-\widetilde{\mathbf U})\partial_z\mathbf U\cdot(\nabla_x\mathbf U-\nabla_x\mathbf u)dxdz
+\int_\Omega\chi_{\rho\geq r}\rho\widetilde{\mathbf U}\partial_z\mathbf U\cdot(\nabla_x\mathbf U-\nabla_x\mathbf u)dxdz\nonumber\\
&=K_1+K_2,\label{333}
\end{align}
where
\begin{align}
K_2&\leq \int_\Omega\chi_{\rho\geq r}\rho^{\frac{\gamma}{2}}\widetilde{\mathbf U}\partial_z\mathbf U\cdot(\nabla_x\mathbf U-\nabla_x\mathbf u)dxdz\nonumber\\
&\leq \|\chi_{\rho\geq r}\rho^{\frac{\gamma}{2}}\|_{L^2(\Omega)}
\|\widetilde{\mathbf U}\partial_z\mathbf U\|_{L^\infty(\Omega)}
\|\nabla_x\mathbf U-\nabla_x\mathbf u\|_{L^2(\Omega)}\nonumber\\
&\leq C\|\chi_{\rho\geq r}\rho^{\frac{\gamma}{2}}\|^2_{L^2(\Omega)}
+\delta\|\nabla_x\mathbf U-\nabla_x\mathbf u\|^2_{L^2(\Omega)}\nonumber\\
&\leq C\mathcal{E}(\rho,\mathbf u|r,\mathbf U)(t)+\delta\|\nabla_x\mathbf U-\nabla_x\mathbf u\|^2_{L^2(\Omega)}.
\end{align}

Next, by virtue of H\"{o}lder inequality, we get
\begin{align*}
K_1&\leq \|\chi_{\rho\geq r}\rho\|_{L^\gamma(\Omega)}
\|\chi_{\rho\geq r}(\widetilde{\mathbf u}-\widetilde{\mathbf U})\|_{L^3(\Omega)}
\|\partial_z\mathbf U\|_{L^\frac{6\gamma}{\gamma-6}(\Omega)}
\|\nabla_x\mathbf U-\nabla_x\mathbf u\|_{L^2(\Omega)}\nonumber\\
&\leq C\|\chi_{\rho\geq r}\rho\|^2_{L^\gamma(\Omega)}
\|\chi_{\rho\geq r}(\widetilde{\mathbf u}-\widetilde{\mathbf U})\|^2_{L^3(\Omega)}
+\delta\|\nabla_x\mathbf U-\nabla_x\mathbf u\|^2_{L^2(\Omega)}\nonumber\\
&\leq  C\|\chi_{\rho\geq r}\rho\|^2_{L^\gamma(\Omega)}
\|\chi_{\rho\geq r}(\widetilde{\mathbf u}-\widetilde{\mathbf U})\|_{L^2(\Omega)}
\|\chi_{\rho\geq r}(\widetilde{\mathbf u}-\widetilde{\mathbf U})\|_{H^1(\Omega)}
+\delta\|\nabla_x\mathbf U-\nabla_x\mathbf u\|^2_{L^2(\Omega)}\nonumber\\
&\leq C\|\chi_{\rho\geq r}\rho\|^4_{L^\gamma(\Omega)}
\|\chi_{\rho\geq r}(\widetilde{\mathbf u}-\widetilde{\mathbf U})\|^2_{L^2(\Omega)}
+\delta\|\chi_{\rho\geq r}(\widetilde{\mathbf u}-\widetilde{\mathbf U})\|^2_{L^2(\Omega)}
+\delta\|\nabla_x\widetilde{\mathbf U}-\nabla_x\widetilde{\mathbf u}\|^2_{L^2(\Omega)}\nonumber\\
&\hspace{10pt}+\delta\|\partial_z\widetilde{\mathbf U}-\partial_z\widetilde{\mathbf u}\|^2_{L^2(\Omega)}
+\delta\|\nabla_x\mathbf U-\nabla_x\mathbf u\|^2_{L^2(\Omega)},
\end{align*}
where we have used the interpolation inequality
\begin{align*}
\|f\|_{L^3}\leq\|f\|^{\frac{1}{2}}_{L^2}\|f\|^{\frac{1}{2}}_{H^1}.
\end{align*}
According \eqref{a3} and \eqref{3aaa}, we have
\begin{align*}
\|\chi_{\rho\geq r}\rho\|^4_{L^\gamma(\Omega)}
=(\int_{\rho\geq r}\rho^\gamma dxdz)^{\frac{4}{\gamma}}
\leq\mathcal{E}(\rho,\mathbf u|r,\mathbf U)^{\frac{4}{\gamma}}(t),
\end{align*}
and
\begin{align*}
\|\chi_{\rho\geq r}(\widetilde{\mathbf u}-\widetilde{\mathbf U})\|^2_{L^2(\Omega)}
=\int_{\rho\geq r}|\widetilde{\mathbf u}-\widetilde{\mathbf U}|^2dxdz
=\int_{\rho\geq r}\frac{1}{\rho}\rho|\widetilde{\mathbf u}-\widetilde{\mathbf U}|^2dxdz
\leq \frac{1}{\|r\|_{\infty(\Omega)}}\mathcal{E}(\rho,\mathbf u|r,\mathbf U)(t).
\end{align*}

Similar to the estimate of \eqref{3aaa}, we obtain
\begin{align*}
\|\nabla_x\widetilde{\mathbf U}-\nabla_x\widetilde{\mathbf u}\|^2_{L^2(\Omega)}
\leq \|\nabla_x\mathbf U-\nabla_x\mathbf u\|^2_{L^2(\Omega)},
\hspace{5pt}\|\partial_z\widetilde{\mathbf U}-\partial_z\widetilde{\mathbf u}\|^2_{L^2(\Omega)}
\leq\|\partial_z\mathbf U-\partial_z\mathbf u\|^2_{L^2(\Omega)}.
\end{align*}

%where we have used the fact the time $T$ is local and finite.

Combining the above estimates, we get
\begin{align*}
\int^\tau_0K_1dt\leq C\int^\tau_0h(t)\mathcal{E}(\rho,\mathbf u|r,\mathbf U)(t)dt
+\delta\int^\tau_0\|\nabla_x\mathbf U-\nabla_x\mathbf u\|^2_{L^2(\Omega)}+
\|\partial_z\mathbf U-\partial_z\mathbf u\|^2_{L^2(\Omega)}dt,
\end{align*}
where $h(t)\in L^1(0,T)$.

Using the same method we estimate the remaining terms. Therefore, we conclude that
\begin{align*}
\mathcal{E}&(\rho,\mathbf{u}|r,\mathbf U)|^{t=\tau}_{t=0}+\int^\tau_0\int_\Omega\big{(}\mu\nabla_x\mathbf u\cdot(\nabla_x\mathbf u-\nabla_x\mathbf U)+\lambda\partial_z\mathbf u(\partial_z\mathbf u-\partial_z\mathbf U)\big{)}dxdzdt\nonumber\\
&\leq C\int^\tau_0h(t)\mathcal{E}(\rho,\mathbf{u}|r,\mathbf U)dt+\delta\int^\tau_0\|\nabla_x\mathbf U-\nabla_x\mathbf u\|^2_{L^{2}(\Omega)}
+\|\partial_z\mathbf U-\partial_z\mathbf u\|^2_{L^{2}(\Omega)}dt\nonumber\\
&\hspace{8pt}+\int^\tau_0\int_\Omega\frac{\rho}{r}(\mathbf U-\mathbf u)(\mu\Delta_x\mathbf U+\lambda\partial_{zz}\mathbf U)dxdzdt
-\int^\tau_0\int_\Omega\frac{\rho}{r}(\mathbf U-\mathbf u)\nabla_xp(r)dxdzdt\nonumber\\
&\hspace{8pt}-\int^\tau_0\int_{\Omega}P''(r)(\rho\partial_tr+\rho\mathbf u\nabla_xr)dxdzdt
+\int^\tau_0\int_\Omega\partial_tp(r)dxdzdt-\int^\tau_0\int_\Omega p(\rho)\text{div}_x\mathbf Udxdzdt.
\end{align*}

Then we deduce that
\begin{align}
\mathcal{E}&(\rho,\mathbf{u}|r,\mathbf U)|^{t=\tau}_{t=0}+\int^\tau_0\int_\Omega\big{(}\mu(\nabla_x\mathbf u-\nabla_x\mathbf U):(\nabla_x\mathbf u-\nabla_x\mathbf U)+\lambda(\partial_z\mathbf u-\partial_z\mathbf U)^2\big{)}dxdzdt\nonumber\\
&\leq C\int^\tau_0h(t)\mathcal{E}(\rho,\mathbf{u}|r,\mathbf U)dt+\delta\int^\tau_0\|\nabla_x\mathbf U-\nabla_x\mathbf u\|^2_{L^{2}(\Omega)}
+\|\partial_z\mathbf U-\partial_z\mathbf u\|^2_{L^{2}(\Omega)}dt\nonumber\\
&\hspace{8pt}+\int^\tau_0\int_\Omega(\frac{\rho}{r}-1)(\mathbf U-\mathbf u)(\mu\Delta_x\mathbf U+\lambda\partial_{zz}\mathbf U)dxdzdt
-\int^\tau_0\int_\Omega\frac{\rho}{r}(\mathbf U-\mathbf u)\nabla_xp(r)dxdzdt\nonumber\\
&\hspace{8pt}-\int^\tau_0\int_{\Omega}P''(r)(\rho\partial_tr+\rho\mathbf u\nabla_xr)dxdzdt
+\int^\tau_0\int_\Omega\partial_tp(r)dxdzdt-\int^\tau_0\int_\Omega p(\rho)\text{div}_x\mathbf Udxdzdt.
\end{align}

It is easy to check that
\begin{align}
-\int^\tau_0&\int_\Omega\frac{\rho}{r}(\mathbf U-\mathbf u)\nabla_xp(r)+p(\rho)\text{div}_x\mathbf U+P''(r)(\rho\partial_tr+\rho\mathbf u\nabla_xr)dxdzdt
+\int^\tau_0\int_\Omega\partial_tp(r)dxdzdt\nonumber\\
&=-\int^\tau_0\int_\Omega(\rho-r)P''(r)\partial_tr+P''(r)\rho\mathbf u\cdot\nabla_xr+\rho P''(r)(\mathbf U-\mathbf u)\cdot\nabla_xr+p(\rho)\text{div}_x\mathbf Udxdzdt\nonumber\\
&=-\int^\tau_0\int_\Omega(\rho-r)P''(r)\partial_tr+P''(r)\rho\mathbf U\cdot\nabla_xr+p(\rho)\text{div}_x\mathbf Udxdzdt\nonumber\\
&=-\int^\tau_0\int_\Omega\rho P''(r)(\partial_tr+\mathbf U\cdot\nabla_xr)
-rP''(r)\partial_tr+p(\rho)\text{div}_x\mathbf Udxdzdt\nonumber\\
&=-\int^\tau_0\int_\Omega\rho P''(r)(-r\text{div}_x\mathbf U-r\partial_zW)
-rP''(r)\partial_tr+p(\rho)\text{div}_x\mathbf Udxdzdt\nonumber\\
&=-\int^\tau_0\int_\Omega\text{div}_x\mathbf U\big{(}p(\rho)-p'(r)(\rho-r)-p(r)\big{)}dxdzdt
+\int^\tau_0\int_\Omega p'(r)(\rho-r)\partial_zWdxdzdt,
\end{align}
where we have used the fact that $\partial_tr+\text{div}_x\mathbf Ur+\mathbf U\cdot\nabla_xr+r\partial_zW=0$.

Recalling the boundary condition $W|_{z=0,1}=0$, we have
\begin{align}
\int^\tau_0\int_\Omega p'(r)(\rho-r)\partial_zWdxdzdt
=\int^\tau_0dt\int_{\mathbb{T}^2}(\int^1_0\partial_zWdz)p'(r)(\rho-r)dx=0.
\end{align}

Moreover, we can use the method as \cite{kr} Section 6.3 to get
\begin{align}
\int_\Omega&(\frac{\rho}{r}-1)(\mathbf U-\mathbf u)(\mu\Delta_x\mathbf U+\lambda\partial_{zz}\mathbf U)dxdz\nonumber\\
&\leq C\mathcal{E}(\rho,\mathbf u|r,\mathbf U)+\delta\|\nabla_x\mathbf u-\nabla_x\mathbf U\|^2_{L^2}
+\delta\|\partial_z\mathbf u-\partial_z\mathbf U\|^2_{L^2}.
\end{align}

Putting $(3.10)-(3.13)$ together, we have
\begin{align}
\mathcal{E}(\rho,\mathbf u|r,\mathbf U)(\tau)\leq C\int^\tau_0h(t)\mathcal{E}(\rho,\mathbf u|r,\mathbf U)(t)dt.
\end{align}

Then applying the Gronwall's inequality, we finish the proof of Theorem 2.1.

\vskip 0.5cm
%%%%%%%%%%%%%%%%%%%%%%%%%%%%%%%%%%%%%%%%%%%%%%%%%%%%%%%%%%%%%%%%
%%%%%%%%%%%%%%%%%%%%%%%%%%%%%%%%%%%%%%%%%%%%%%%%%%%%%%%%%%%%%%%%
\noindent {\bf Part II: Singular limit of CPE}
\vskip 0.2cm

This part is devoted to studying the singular limit of the CPE in the case of well-prepared initial data.

\section{Preliminaries and main result}
From the notable survey paper by Klein, see \cite{Klein}, singular limits of fluids play an important role in mathematics, physics and meteorology. We consider the following scale CPE system with Coriolis forces:
\begin{eqnarray}
\left\{
\begin{array}{llll}  \partial_{t}\rho_\epsilon+\text{div}_x(\rho_\epsilon \mathbf{u}_\epsilon)+\partial_z(\rho_\epsilon w_\epsilon)=0, \\
\partial_t(\rho_\epsilon\mathbf{u}_\epsilon)+\textrm{div}_x(\rho_\epsilon\mathbf{u}_\epsilon\otimes\mathbf{u}_\epsilon)
+\partial_x(\rho_\epsilon\mathbf u_\epsilon w_\epsilon)
+\rho_\epsilon\mathbf u_\epsilon\times\omega+\frac{1}{\epsilon^2}\nabla_x p(\rho_\epsilon)=\mu\Delta_x\mathbf u_\epsilon+\lambda\partial^2_{zz}\mathbf u_\epsilon,\\
\partial_zp(\rho_\epsilon)=0,
\end{array}\right.\label{4a}
\end{eqnarray}
where $\epsilon$ represents the Mach number, $\omega=(0,0,1)$ is the rotation axis. The boundary conditions and pressure are the same as (1.2) and (1.4). Problem \eqref{4a} is supplemented with initial data
\begin{align}
\rho_\epsilon (0, \cdot) = \rho_{0, \epsilon} =\overline{\rho} + \epsilon \rho^{(1)}_{0, \epsilon},\hspace{5pt}
\mathbf{u}_\epsilon (0, \cdot) = \mathbf{u}_{0, \epsilon},
\end{align}
where the constant $\overline{\rho}$ in (4.2) can be taken arbitrary.

There is a quite broad consensus that the compressible flows become incompressible in the low Mach number limit. In the following sections, we assume $\rho=\rho_\epsilon$ and $\mathbf u=\mathbf u_\epsilon$. In this part, our goal is to study system \eqref{4a} in the case of singular limit $\epsilon\rightarrow0$, meaning the inviscid, incompressible limit. Precisely speaking, we want to show that the weak solutions of CPE converge to the incompressible PE system.

\subsection{Target equation}

The expected limit problem reads
\begin{align}
&\partial_t\mathbf{V}+(\mathbf{V}\cdot\nabla_x)\mathbf{V}+\partial_z\mathbf VW+\mathbf V^{\perp}+\nabla_x\Pi=0,\nonumber\\
&\text{div}_x\mathbf{V}+\partial_zW=0,\nonumber\\
&\partial_z\Pi=0,\label{4b}
\end{align}
where $\mathbf V^{\perp}=(v_2,-v_1)$ and the $\Pi$ is the pressure. We supplement the system with the initial condition
\begin{align*}
\mathbf V|_{t=0}=\mathbf V_0.
\end{align*}

As shown by Kukavica et al. \cite{kukavica}, the problem \eqref{4b} possesses a local unique analytic solution $\mathbf V$ and $\Pi$ for some $T>0$ and any initial solution
\begin{align}
\mathbf V_0\in C^{\infty}(\Omega),\hspace{3pt}\int^1_0\text{div}\mathbf V_0dz=0.
\end{align}

\subsection{Relative energy inequality}
According to the previous definition, we define the relative entropy functional,
\begin{align}
\mathcal{E}(\rho,\mathbf{u}|r, \mathbf{V})=\int_{\Omega}[\frac{1}{2}\rho|\mathbf u-\mathbf V|^2+\frac{1}{\epsilon^2}(P(\rho)-P'(r)(\rho-r)-P(r))]dxdz,
\end{align}
where $r$ and $\mathbf V$ are continuously differentiable, it is something not understandable "text functions". The following relation can be deduced
\begin{align}
\mathcal{E}&(\rho,\mathbf{u}|r,\mathbf V)|^{t=\tau}_{t=0}+\int^\tau_0\int_\Omega\big{(}\mu\nabla_x\mathbf u\cdot(\nabla_x\mathbf u-\nabla_x\mathbf V)+\lambda\partial_z\mathbf u(\partial_z\mathbf u-\partial_z\mathbf V)\big{)}dxdzdt\nonumber\\
&\leq\int^\tau_0\int_{\Omega}\rho(\mathbf V-\mathbf u)\partial_t\mathbf V+\rho\mathbf u(\mathbf V-\mathbf u)\cdot\nabla_x\mathbf V
+\rho w(\mathbf V-\mathbf u)\partial_z\mathbf V-\frac{1}{\epsilon^2}p(\rho)\text{div}_x\mathbf Vdxdzdt\nonumber\\
&\hspace{15pt}-\frac{1}{\epsilon^2}\int^\tau_0\int_{\Omega}P''(r)(\rho\partial_tr+\rho\mathbf u\nabla_xr)dxdzdt
+\frac{1}{\epsilon^2}\int^\tau_0\int_{\Omega}\partial_tp(r)dxdzdt\nonumber\\
&\hspace{15pt}-\int^\tau_0\int_\Omega(\rho\mathbf u\times\omega)\cdot(\mathbf V-\mathbf u)dxdzdt,\label{4c}
\end{align}
for any $r,\mathbf V$$\in C'([0,T]\times\Omega)$, $r>0$.

%%%%%%%%%%%%%%%%%%%%%%%%%%%%%%%%%%%%%%%%%%%%%
%%%%%%%%%%%%%%%%%%%%%%%%%%%%%%%%%%%%%%%%%%%%%
\subsection{Main result}
The second result concerns the singular limit.

\begin{theorem}
Let $\gamma>6$, and $(\rho,\mathbf u,w)$ be a weak solution of the scaled system \eqref{4a} on a time interval $(0,T)$ with well-prepared initial data satisfying the following assumptions
\begin{align}
&\|\rho^{(1)}_{0,\epsilon}\|_{L^\infty(\Omega)}+\|\mathbf u_{0,\epsilon}\|_{L^\infty(\Omega)}\leq D,\nonumber\\
&\frac{\rho_{0,\epsilon}-\overline{\rho}}{\epsilon}\rightarrow 0\hspace{3pt}\text{in}\hspace{3pt}L^{1}(\Omega),\hspace{8pt}
\mathbf{u}_{0,\epsilon}\rightarrow \mathbf{V}_0\hspace{3pt}\text{in}\hspace{3pt}L^{2}(\Omega).
\end{align}

Let $\mathbf V$ be the unique analytic solution of the target problem \eqref{4b}. Suppose that $T<T_{\max}$, where $T_{\max}$ denotes the maximal life-span of the regular solution to the incompressible PE system \eqref{4b} with initial data $\mathbf V_0$, then
\begin{align}
\sup_{t\in[0,T]}&\int_\Omega[\rho|\mathbf u-\mathbf V|^2+\frac{1}{\epsilon^2}(P(\rho)-P'(\overline{\rho})(\rho-\overline{\rho})-P(\overline{\rho}))]\nonumber\\
&\leq C[\epsilon+\mu+\lambda+\int_\Omega|\mathbf u_{0,\epsilon}-\mathbf V_0|^2],
\end{align}
where the constant $C$ depends on the initial data $\rho_{0}$, $\mathbf u_{0}$, $\mathbf V_0$ and $T$, and the size $D$ of the initial data perturbation. The constant $\overline{\rho}$ can be taken arbitrary.
\end{theorem}

\begin{remark}
Theorem 4.1 yields that $\rho_\epsilon$ and $\mathbf u_\epsilon$ converge to the solution of target system in the regime of $\epsilon\rightarrow0$ and $\mu,\lambda\rightarrow0$ for the well-prepared initial data, in other words, the expression of the right hand of (4.8) tends to zero.
\end{remark}

%%%%%%%%%%%%%%%%%%%%%%%%%%%%%%%%%%%%%%%%%%%%%%%%%%%%%%%%%%%%%%%%
%%%%%%%%%%%%%%%%%%%%%%%%%%%%%%%%%%%%%%%%%%%%%%%%%%%%%%%%%%%%%%%%
\subsection{Uniform bounds}
Before proving Theorems 4.1, we derive uniforms bounds of weak solutions $(\rho,\mathbf u)$. Here and hereafter, the constant $C$ denotes a positive constant, independent on $\epsilon$, that will not have the same value when used in different parts of text. The following uniform bounds are derived from the relative energy inequality \eqref{4c}, if we take $r=\overline{\rho}$ and $\mathbf U=0$:
\begin{align}
&ess\sup_{t\in(0,T)}||\frac{\rho-\overline{\rho}}{\epsilon}||_{L^2(\Omega)+L^\gamma(\Omega)}\leq C,\nonumber\\
&ess\sup_{t\in(0,T)}||\sqrt{\rho}\mathbf u||_{L^2(\Omega)}\leq C,\hspace{5pt}
\sqrt{\mu}||\nabla_x\mathbf u||_{L^2((0,T)\times\Omega)}
+\sqrt{\lambda}||\partial_z\mathbf u||_{L^2((0,T)\times\Omega)}\leq C.
\end{align}

%%%%%%%%%%%%%%%%%%%%%%%%%%%%%%%%%%%%%%%%%%%%%%%%%%
%%%%%%%%%%%%%%%%%%%%%%%%%%%%%%%%%%%%%%%%%%%%%%%%%%%%
\section{Convergence of well-prepared initial data}

The proof of convergence is based on the ansatz
\begin{equation}
r=\overline{\rho},\hspace{5pt} \mathbf{U}=\mathbf{V},
\end{equation}
in the relative energy inequality \eqref{4c}, where $\mathbf V$ is the analytic solution of the target problem \eqref{4b}. The corresponding relative energy inequality reads as:
\begin{align}
\mathcal{E}&(\rho,\mathbf{u}|\overline{\rho}, \mathbf V)(\tau)+\int^\tau_0\int_\Omega\mu(\nabla_x\mathbf u-\nabla_x\mathbf V\big{)}:(\nabla_x\mathbf u -\nabla_x\mathbf V)
+\lambda(\partial_z\mathbf u-\partial_z\mathbf V)^2dxdzdt\nonumber\\
&\leq\mathcal{E}(\rho,\mathbf{u}|\overline{\rho}, \mathbf V)(0)
+\int^\tau_0\int_{\Omega}\rho(\mathbf V-\mathbf u)\partial_t\mathbf V+\rho\mathbf u(\mathbf V-\mathbf u)\cdot\nabla_x\mathbf V+\rho w(\mathbf V-\mathbf u)\partial_z\mathbf Vdxdzdt\nonumber\\
&\hspace{10pt}-\frac{1}{\epsilon^2}\int^\tau_0\int_\Omega p(\rho)\text{div}_x\mathbf Vdxdzdt
-\int^\tau_0\int_\Omega(\rho\mathbf u\times\omega)\cdot(\mathbf V-\mathbf u)dxdzdt\nonumber\\
&\hspace{10pt}+\int^\tau_0\int_\Omega\mu\nabla_x\mathbf V(\nabla_x\mathbf u-\nabla_x\mathbf V)dxdzdt
+\int^\tau_0\int_\Omega\lambda\partial_z\mathbf V(\partial_z\mathbf u-\partial_z\mathbf V)dxdzdt.
\end{align}

First we deal with initial data and viscous term. It is easy to computer the initial relative energy inequality:
\begin{align}
\mathcal{E} (\rho_{},\mathbf{u}_{}|\overline{\rho}, \mathbf V)|_{t=0}\leq C\int_\Omega[|\mathbf u_{0,\epsilon}-\mathbf V_0|^2+|\rho_{0,\epsilon}-\overline{\rho}|^2]dx,
\end{align}
and viscous term
\begin{align}
&\mu\int^\tau_0\int_\Omega\nabla_x\mathbf V(\nabla_x\mathbf u-\nabla_x\mathbf V)dxdzdt\leq
\int^\tau_0\frac{\mu}{2}\|\nabla_x\mathbf u-\nabla_x\mathbf V\|^2_{L^2(\Omega)}
+\frac{\mu}{2}\|\nabla_x\mathbf V\|^2_{L^2(\Omega)}dt,\nonumber\\
&\lambda\int^\tau_0\int_\Omega\partial_z\mathbf V(\partial_z\mathbf u-\partial_z\mathbf V)dxdzdt\leq
\int^\tau_0\frac{\lambda}{2}\|\nabla_x\mathbf u-\nabla_x\mathbf V\|^2_{L^2(\Omega)}
+\frac{\lambda}{2}\|\partial_z\mathbf V\|^2_{L^2(\Omega)}dt.
\end{align}

Next, we consider the remaining terms. Utilizing $(4.3)_1$, we get that
\begin{align}
\int^\tau_0\int_{\Omega}\rho&(\mathbf V-\mathbf u)\partial_t\mathbf V+\rho\mathbf u(\mathbf V-\mathbf u)\cdot\nabla_x\mathbf V+\rho w(\mathbf V-\mathbf u)\partial_z\mathbf Vdxdzdt\nonumber\\
&=\int^\tau_0\int_{\Omega}\rho(\mathbf V-\mathbf u)(\partial_t\mathbf V+(\mathbf V\cdot\nabla_x)\mathbf V+W\partial_z\mathbf V)+\rho(\mathbf u-\mathbf V)(\mathbf V-\mathbf u)\nabla_x\mathbf V\nonumber\\
&\hspace{30pt}+\rho(\mathbf V-\mathbf u)(w-W)\partial_z\mathbf Vdxdzdt\nonumber\\
&=\int^\tau_0\int_{\Omega}\rho(\mathbf u-\mathbf V)(\nabla_x\Pi+\mathbf V^\bot)+\rho(\mathbf u-\mathbf V)(\mathbf V-\mathbf u)\nabla_x\mathbf V
+\rho(\mathbf V-\mathbf u)(w-W)\partial_z\mathbf Vdxdzdt.
\end{align}

It is easy to check that
\begin{align*}
\int^\tau_0\int_\Omega\rho(\mathbf u-\mathbf V)(\mathbf V-\mathbf u)\nabla_x\mathbf Vdxdzdt\leq C\int^\tau_0\mathcal{E}(\rho,\mathbf{u}|\overline{\rho}, \mathbf V)(t)dt.
\end{align*}

Next, we estimate the term $\int^\tau_0\int_{\Omega}\rho\mathbf V\cdot\nabla_x\Pi dxdzdt$, and rewrite in the form
\begin{align}
\int^\tau_0\int_{\Omega}\rho\mathbf V\cdot\nabla_x\Pi dxdzdt=
\epsilon\int^\tau_0\int_{\Omega}\frac{\rho-\overline{\rho}}{\epsilon}\mathbf V\cdot\nabla_x\Pi dxdzdt+
\overline{\rho}\int^\tau_0\int_{\Omega}\mathbf V\cdot\nabla_x\Pi dxdzdt,
\end{align}

The second term on the right side of (5.6) is estimated as:
\begin{align*}
\int^\tau_0\int_{\Omega}\mathbf V\nabla_x\Pi dxdzdt
=-\int^\tau_0\int_{\Omega}\text{div}_x\mathbf V\Pi dxdzdt
=\int^\tau_0\int_{\Omega}\partial_zW\Pi dxdzdt=0,
\end{align*}
where we have used the fact that $\Pi$ is independent of $z$. We deduce from the energy inequality that
\begin{align}
\int_{\Omega}\frac{1}{\epsilon^2}(P(\rho)-P'(r)(\rho-r)-P(r))dxdz\leq C, \hspace{5pt}\text{uniformly as}
\hspace{3pt}\epsilon\rightarrow0.
\end{align}

Similar to the previous analysis, it is enough to establish a uniform bound
\begin{align*}
\int_\Omega\frac{\rho-\overline{\rho}}{\epsilon}dxdz\leq C.
\end{align*}

As we know that the pressure $\Pi$ is analytic, so that the rightmost integral of (5.6) can be vanished as $\epsilon\rightarrow0$.

From the previous definition of dissipative weak solutions, we choose $\Pi$ as the test function, so that
\begin{align*}
\int^\tau_0\int_{\Omega}&\rho\mathbf u\nabla_x\Pi dxdzdt\\
&=[\int_\Omega\rho\Pi dxdz]|^{t=\tau}_{t=0}-\int^\tau_0\int_\Omega\rho\partial_t\Pi dxdzdt
-\int^\tau_0\int_\Omega\rho w\partial_z\Pi dxdzdt\\
&=\epsilon[\int_\Omega\frac{\rho-\overline{\rho}}{\epsilon}\Pi dxdz]|^{t=\tau}_{t=0}
-\epsilon\int^\tau_0\int_\Omega\frac{\rho-\overline{\rho}}{\epsilon}\partial_t\Pi dxdzdt\rightarrow0,
\hspace{3pt}\text{as}\hspace{3pt}\epsilon\rightarrow0.
\end{align*}

Compared with Navier-Stokes equations, the pressure term in PE system is easy to estimate. By virtue of incompressible condition and $(4.3)_3$, we have that
\begin{align*}
-\frac{1}{\epsilon^2}\int^\tau_0\int_\Omega p(\rho)\text{div}_x\mathbf Vdxdzdt=\frac{1}{\epsilon^2}\int^\tau_0\int_\Omega p(\rho)\partial_zWdxdzdt=0.
\end{align*}

Moreover, we find that
\begin{align*}
\int_{\Omega}\rho(\mathbf u-\mathbf V)\cdot\mathbf V^\bot dxdz+\int_\Omega(\rho\mathbf u\times\omega)\cdot(\mathbf V-\mathbf u)dxdz=0.
\end{align*}

Now, utilizing \eqref{b1}, we deal the complex nonlinear term
\begin{align*}
\int^\tau_0\int_{\Omega}&\rho(\mathbf V-\mathbf u)(w-W)\partial_z\mathbf V dxdzdt\\
&=\int^\tau_0\int_{\Omega}(\mathbf V-\mathbf u)\partial_z\mathbf V\big{(}
-\text{div}_x(\rho\widetilde{\mathbf u})+z\text{div}_x(\rho\overline{\mathbf u})
-\rho W \big{)}dxdzdt.
\end{align*}

These nonlinear terms are estimated one by one
\begin{align}
-\int_{\Omega}(\mathbf V-\mathbf u)\partial_z\mathbf V
\text{div}_x(\rho\widetilde{\mathbf u})dxdz
=\int_{\Omega}\rho\widetilde{\mathbf u}(\nabla_x\mathbf V-\nabla_x\mathbf u)\cdot\partial_z\mathbf Vdxdz
+\int_{\Omega}\rho\widetilde{\mathbf u}(\mathbf V-\mathbf u)\cdot\partial_z\nabla_x\mathbf Vdxdz.\label{55a}
\end{align}

From the incompressible condition, it follows that $W=-\int^z_0\text{div}_x\mathbf V(x,s,t)ds$. We define $\widetilde{\mathbf V}=\int^z_0\mathbf V(x,s,t)ds$ and get
\begin{align}
\int_{\Omega}&\rho\widetilde{\mathbf u}(\nabla_x\mathbf V-\nabla_x\mathbf u)\cdot\partial_z\mathbf Vdxdz\nonumber\\
&=\int_{\Omega}\chi_{\rho\leq\overline{\rho}}\rho\widetilde{\mathbf u}(\nabla_x\mathbf V-\nabla_x\mathbf u)\partial_z\mathbf Vdxdz
+\int_{\Omega}\chi_{\rho\geq\overline{\rho}}\rho\widetilde{\mathbf u}(\nabla_x\mathbf V-\nabla_x\mathbf u)\partial_z\mathbf Vdxdz\nonumber\\
&=\int_{\Omega}\chi_{\rho\leq\overline{\rho}}\rho(\widetilde{\mathbf u}-\widetilde{\mathbf V})(\nabla_x\mathbf V-\nabla_x\mathbf u)\partial_z\mathbf Vdxdz
+\int_{\Omega}\chi_{\rho\leq\overline{\rho}}\rho\widetilde{\mathbf V}(\nabla_x\mathbf V-\nabla_x\mathbf u)\partial_z\mathbf Vdxdz\nonumber\\
&\hspace{5pt}+\int_{\Omega}\chi_{\rho\geq\overline{\rho}}\rho\widetilde{\mathbf u}(\nabla_x\mathbf V-\nabla_x\mathbf u)\partial_z\mathbf Vdxdz\label{5a}
\end{align}

The foremost two terms on the right side of \eqref{5a} can be handed as \eqref{2a}
\begin{align}
\int_{\Omega}&\chi_{\rho\leq\overline{\rho}}\rho(\widetilde{\mathbf u}-\widetilde{\mathbf V})(\nabla_x\mathbf V-\nabla_x\mathbf u)\partial_z\mathbf Vdxdz
+\int_{\Omega}\chi_{\rho\leq\overline{\rho}}\rho\widetilde{\mathbf V}(\nabla_x\mathbf V-\nabla_x\mathbf u)\partial_z\mathbf Vdxdz\nonumber\\
&=\int_{\Omega}\chi_{\rho\leq\overline{\rho}}\rho(\widetilde{\mathbf u}-\widetilde{\mathbf V})(\nabla_x\mathbf V-\nabla_x\mathbf u)\partial_z\mathbf Vdxdz
+\int_{\Omega}\chi_{\frac{\overline{\rho}}{2}<\rho\leq\overline{\rho}}\rho\widetilde{\mathbf V}(\nabla_x\mathbf V-\nabla_x\mathbf u)\partial_z\mathbf Vdxdz\nonumber\\
&\hspace{5pt}+\int_{\Omega}\chi_{\rho\leq\frac{\overline{\rho}}{2}}\rho\widetilde{\mathbf V}(\nabla_x\mathbf V-\nabla_x\mathbf u)\partial_z\mathbf Vdxdz\nonumber\\
&\leq\delta\|\nabla_x\mathbf V-\nabla_x\mathbf u\|_{L^2(\Omega)}^2
+C\mathcal{E}(\rho,\mathbf{u}|\overline{\rho}, \mathbf V)(t).
\end{align}

On the other hand, following \eqref{333}, we have
\begin{align}
\int^\tau_0\int_{\Omega}&\chi_{\rho\geq\overline{\rho}}\rho\widetilde{\mathbf u}(\nabla_x\mathbf V-\nabla_x\mathbf u)\partial_z\mathbf Vdxdz\nonumber\\
&=\int^\tau_0\int_{\Omega}\chi_{\rho\geq\overline{\rho}}\rho(\widetilde{\mathbf u}-\widetilde{\mathbf V})(\nabla_x\mathbf V-\nabla_x\mathbf u)\partial_z\mathbf Vdxdz
+\int^\tau_0\int_{\Omega}\chi_{\rho\geq\overline{\rho}}\rho\widetilde{\mathbf V}(\nabla_x\mathbf V-\nabla_x\mathbf u)\partial_z\mathbf Vdxdz\nonumber\\
&\leq C\int^\tau_0h(t)\mathcal{E}(\rho,\mathbf u|r,\mathbf U)(t)dt
+\delta\int^\tau_0\|\nabla_x\mathbf V-\nabla_x\mathbf u\|^2_{L^2(\Omega)}
+\|\partial_z\mathbf V-\partial_z\mathbf u\|^2_{L^2(\Omega)}dt.
\end{align}

Similarly, the second nonlinear term on the right side of \eqref{55a} is divided into two parts:
\begin{align}
\int_{\Omega}\rho\widetilde{\mathbf u}(\mathbf V-\mathbf u)\cdot\partial_z\nabla_x\mathbf Vdxdz
&=\int_{\Omega}\rho(\widetilde{\mathbf u}-\widetilde{\mathbf V})(\mathbf V-\mathbf u)\partial_z\nabla_x\mathbf Vdxdz
+\int_{\Omega}\rho\widetilde{\mathbf V}(\mathbf V-\mathbf u)\partial_z\nabla_x\mathbf Vdxdz.\nonumber
\end{align}

Utilizing the similar estimates in (3.6), we have
\begin{align}
\int_{\Omega}\rho(\widetilde{\mathbf u}-\widetilde{\mathbf V})(\mathbf V-\mathbf u)\partial_z\nabla_x\mathbf Vdxdz
\leq C\mathcal{E}(\rho,\mathbf{u}|\overline{\rho}, \mathbf V)(t).
\end{align}

Moreover, similar to \eqref{33c}, we get
\begin{align}
\int_{\Omega}&\rho\widetilde{\mathbf V}(\mathbf V-\mathbf u)\partial_z\nabla_x\mathbf Vdxdz\nonumber\\
&=\int_{\Omega}\chi_{\rho\leq\frac{\overline{\rho}}{2}}\rho\widetilde{\mathbf V}(\mathbf V-\mathbf u)\partial_z\nabla_x\mathbf Vdxdz
+\int_{\Omega}\chi_{\frac{\overline{\rho}}{2}<\rho<\overline{\rho}}\rho\widetilde{\mathbf V}(\mathbf V-\mathbf u)\partial_z\nabla_x\mathbf Vdxdz\nonumber\\
&\hspace{5pt}+\int_{\Omega}\chi_{\rho\geq\overline{\rho}}\rho\widetilde{\mathbf V}(\mathbf V-\mathbf u)\partial_z\nabla_x\mathbf Vdxdz\nonumber\\
&\leq \delta\|\nabla_x\mathbf V-\nabla_x\mathbf u\|_{L^2(\Omega)}^2
+\delta\|\partial_z\mathbf V-\partial_z\mathbf u\|_{L^2(\Omega)}^2
+C\mathcal{E}(\rho,\mathbf{u}|\overline{\rho}, \mathbf V)(t),
\end{align}
and
\begin{align}
\int_{\Omega}(\mathbf V-\mathbf u)\partial_z\mathbf Vz\text{div}_x(\rho\overline{\mathbf u})dxdz
\leq \delta\|\nabla_x\mathbf V-\nabla_x\mathbf u\|_{L^2(\Omega)}^2
+\delta\|\partial_z\mathbf V-\partial_z\mathbf u\|_{L^2(\Omega)}^2
+C\mathcal{E}(\rho,\mathbf{u}|\overline{\rho}, \mathbf V)(t).
\end{align}

The last term can be estimated as
\begin{align*}
\int_{\Omega}&\rho(\mathbf V-\mathbf u)\partial_z\mathbf VWdxdz\\
&=\int_\Omega\chi_{\rho\leq\frac{\overline{\rho}}{2}}\rho(\mathbf V-\mathbf u)\partial_z\mathbf VWdxdz+
\int_\Omega\chi_{\rho\geq\overline{\rho}}\rho(\mathbf V-\mathbf u)\partial_z\mathbf VWdxdz\\
&\hspace{8pt}+\int_\Omega\chi_{\frac{\overline{\rho}}{2}<\rho<\overline{\rho}}\rho(\mathbf V-\mathbf u)\partial_z\mathbf VWdxdz\\
& \leq\delta\|\nabla_x\mathbf V-\nabla_x\mathbf u\|_{L^2(\Omega)}^2
+\delta\|\partial_z\mathbf V-\partial_z\mathbf u\|_{L^2(\Omega)}^2
+C\mathcal{E}(\rho,\mathbf{u}|\overline{\rho}, \mathbf V)(t).
\end{align*}

%%%%%%%%%%%%%%%%%%%%%%%%%%%%%%%%%%%%%%%%%%%%%%%%%%%%%%

Combining the above estimates together and using Grownwall inequality, we prove Theorem 4.1.

%%%%%%%%%%%%%%%%%%%%%%%%%%%%%%%%%%%%%%%%%%%%%%%%%%%%%%%%%%%%%%%%%%%%%%%%%%%%%%%%%%%%

\vskip 0.5cm
\noindent {\bf Acknowledgements}

\vskip 0.1cm
We are very much indebted to an anonymous referee for many helpful suggestions. The research of H. G is partially supported by the NSFC Grant  No. 11531006. The research of \v S.N. leading to these results has received funding from the Czech Sciences Foundation (GA\v CR),   GA19-04243S and RVO 67985840.  The research of T.T. is supported by the NSFC Grant No. 11801138. The paper was written when Tong Tang was visiting the Institute of Mathematics of the Czech Academy of Sciences which {hospitality and support} is gladly acknowledged.

%%%%%%%%%%%%%%%%%%%%%%%%%%%%%%%%%%%%%%%%%%%%%%%%%%%%%%%%%


\begin{thebibliography}{99}

\bibitem{b} Y. Brenier, Homogeneous hydrostatic flows with convex velocity profiles, {\em Nonlinearity}, 12 (1999), 495-512.

\bibitem{by} Y. Brenier, Remarks on the derivation of the hydrostatic Euler equations, {\em Bull. Sci. Math.}, 127 (2003), 585-595.

\bibitem{bg} D. Bresch, F. Guill\'{e}n-Gonz\'{a}lez, N. Masmoudi and M. A. Rodr\'{\i}guez-Bellido, On the uniqueness of weak solutions of the two-dimensional primitive equations, {\em Differential Integral Equations}, 16 (2003), 77-94.


\bibitem{b1} D. Bresch, A. Kazhikhov and J. Lemoine, On the two-dimensional hydrostatic Navier-Stokes equations, {\em SIAM J. Math. Anal.}, 36 (2004/05), 796-814.

\bibitem{b2} D. Bresch and B. Desjardins, On the construction of approximate solutions for the 2D viscous shallow water model and for compressible Navier-Stokes models, {\em J. Math. Pures Appl.}, 86 (2006), 362-368.

\bibitem{b3} D. Bresch and B. Desjardins, Existence of global weak solutions to the Navier-Stokes equations for viscous compressible and heat conducting fluids, {\em J. Math. Pures Appl.}, 87 (2007), 57-90.

\bibitem{br} D. Bresch and P. E. Jabin, Global existence of weak solutions for compressible Navier-Stokes equations: thermodynamically unstable pressure and anisotropic viscous stress tensor, {\em Ann. of Math.}, 188 (2018), 577-684.

\bibitem{c1} C. S. Cao and E. S. Titi, Global well-posedness of the three-dimensional viscous primitive equations of large scale ocean and atmosphere dynamics, {\em Ann. of Math.}, 166 (2007), 245-267.

\bibitem{c2} C. S. Cao, J. K. Li and E. S. Titi, Local and global well-posedness of strong solutions to the 3D primitive equations with vertical eddy diffusivity, {\em Arch. Ration. Mech. Anal.}, 214 (2014), 35-76.

\bibitem{c3} C. S. Cao, S. Ibrahim, K. Nakanishi and E. S. Titi, Finite-time blowup for the inviscid primitive equations of oceanic and atmospheric dynamics, {\em Comm. Math. Phys.} 337 (2015), 473-482.

\bibitem{c4} C. S. Cao, J. K. Li and E. S. Titi, Global well-posedness of the three-dimensional primitive equations with only horizontal viscosity and diffusion, {\em Comm. Pure Appl. Math.}, 69 (2016), 1492-1531.

\bibitem{c5} C. S. Cao, J. K. Li and E. S. Titi, Strong solutions to the 3D primitive equations with only horizontal dissipation: near $H^1$ initial data, {\em J. Funct. Anal.}, 272 (2017), 4606-4641.

\bibitem{ch} J. Y. Chemin, B. Desjardins, I. Gallagher and E. Grenier, Mathematical Geophysics: An Introduction
to Rotating Fluids and the Navier-Stokes Equations, Oxford University Press, Oxford, 2006.

\bibitem{D}
C. M. Dafermos,
 The second law of thermodynamics and stability,
{\em Arch. Rational Mech. Anal.}, 70 (1979) 167--179.

\bibitem{d1} A. Debussche, N. Glatt-Holtz, and R. Temam, Local martingale and
pathwise solutions for an abstract fluids model, \emph{Physica D},  240 (2011), 1123-1144.


\bibitem{d2} A. Debussche, N. Glatt-Holtz, R. Temam and M. Ziane, Global existence and regularity for the 3D stochastic primitive equations of the ocean and atmosphere with multiplicative white noise, \emph {Nonlinearity},  316 (2012), 723-76.

\bibitem{dong} Z. Dong, J. Zhai and R. Zhang,  Large deviation principles for 3D stochastic primitive equations, {\em J. Differential Equations}, 263 (2017), 3110-3146.


\bibitem{er1} M. Ersoy, T. Ngom and M. Sy, Compressible primitive equations: formal derivation and stability of weak solutions, {\em Nonlinearity}, 24 (2011), 79-96.


\bibitem{er2} M. Ersoy and T. Ngom,  Existence of a global weak solution to one model of compressible primitive equations, {\em C. R. Math. Acad. Sci. Paris}, 350 (2012), 379-382.


\bibitem{e1} E. Feireisl, Dynamics of viscous compressible fluids, Oxford Lecture Series in Mathematics and its Applications, Oxford University Press, Oxford, 2004.

\bibitem{e2} E. Feireisl, J. B. Jin and A. Novotn\'{y}, Relative entropies, suitable weak solutions, and weak-strong uniqueness for the compressible Navier-Stokes system, {\em J. Math. Fluid Mech.}, 14 (2012), 717-730.

\bibitem{e3} E. Feireisl and A. Novotn\'{y}, Singular limits in thermodynamics of viscous fluids, Advances in Mathematical Fluid Mechanics, Birkh\"{a}user, Basel, 2009.

\bibitem{e} E. Feireisl, I. Gallagher, A. Novotn\'{y}, A singular limit for compressible rotating fluids, {\em SIAM J.
Math. Anal.}, 44 (1) (2012), 192-205.

\bibitem{e4} E. Feireisl, J. B. Jin and A. Novotn\'{y}, Inviscid incompressible limits of strongly stratified fluids, {\em  Asymptot. Anal.}, 89 (2014), 307-329.


\bibitem{gao} H. Gao and C. Sun, Well-posedness of stochastic primitive equations with multiplicative noise in three dimensions, {\em Discrete Contin. Dyn. Syst. Ser. B},  21 (2016), no. 9, 3053-3073.

\bibitem{g} B. V. Gatapov and A. V. Kazhikhov, Existence of a global solution of a model problem of atmospheric dynamics, {\em Siberian Math. J.}, 46 (2005), 805-812.

\bibitem{ge} P. Germain, Weak-strong uniqueness for the isentropic compressible Navier-Stokes system, {\em J. Math. Fluid Mech.}, 13 (2011), 137-146.

\bibitem{gu} F. Guill\'{e}n-Gonz\'{a}lez, N. Masmoudi and M. A. Rodr\'{\i}guez-Bellido, Anisotropic estimates and strong solutions of the primitive equations, {\em Differential Integral Equations}, 14 (2001),  1381-1408.


\bibitem{g1} B. L. Guo and D. W. Huang, Existence of weak solutions and trajectory attractors for the moist atmospheric equations in geophysics, {\em J. Math. Phys.}, 47 (2006), 083508.

\bibitem{g2} B. L. Guo and D. W. Huang, Existence of the universal attractor for the 3-D viscous primitive equations of large-scale moist atmosphere, {\em J. Differential Equations}, 251 (2011), 457-491.

\bibitem{g3} B. L. Guo, D. W. Huang and W. Wang, Diffusion limit of 3D primitive equations of the large-scale ocean under fast oscillating random force, {\em J. Differential Equations}, 259 (2015), 2388-2407.

\bibitem{ju} N. Ju, The global attractor for the solutions to the 3d viscous
primitive equations, {\em Discrete Contin. Dyn. Syst.}, \textbf{17},
(2007), 159-179.

\bibitem{kukavica} I. Kukavica and M. Ziane,  On the regularity of the primitive
equations of the ocean, {\em Nonlinearity}, \textbf{20}, (2007),  2739-2753.


\bibitem{Klein}
R. Klein,
\newblock Scale-dependent models for atmospheric flows,
\newblock In {\em Annual review of fluid mechanics,} Vol. 42, Annu. Rev.
  Fluid. Mech., pages 249--274. Annual Reviews, Palo Alto, CA, 2010.


\bibitem{kr} O. Kreml, \v{S}. Ne\v{c}asov\'{a} and T. Piasecki, Local existence of strong solution and weak-strong uniqueness for the compressible Navier-Stokes system on moving domains, accepted in  {\em Proceedings of the Royal Society of Edinburgh Section A: Mathematics}, DOI: https://doi.org/10.1017/prm.2018.165.



\bibitem{k} I. Kukavica, R. Temam, V. C. Vicol and M. Ziane, Local existence and uniqueness for the hydrostatic Euler equations on a bounded domain, {\em J. Differential Equations}, 250 (2011), 1719-1746.

\bibitem{li} J. Li and Z. P. Xin, Global existence of weak solutions to the barotropic compressible Navier-
Stokes flows with degenerate viscosities, {\em arXiv:1504.06826v2}, 2015.


\bibitem{l1} J. L. Lions, R. Temam and S. H. Wang, On the equations of the large-scale ocean, {\em Nonlinearity}, 5 (1992), 1007-1053.

\bibitem{l2} J. L. Lions, R. Temam and S. H. Wang, New formulations of the primitive equations of atmosphere and applications, {\em Nonlinearity}, 5 (1992), 237-288.

\bibitem{l3} J. L. Lions, R. Temam and S. H. Wang, Mathematical theory for the coupled atmosphere-ocean models, (CAO III), {\em J. Math. Pures Appl.}, (9) 74 (1995), 105-163.

\bibitem{l4} J. L. Lions, O. P. Manley, R. Temam and S. H. Wang, Physical interpretation of the attractor dimension for the primitive equations of atmospheric circulation, {\em J. Atmospheric Sci.}, 54 (1997), 1137-1143.

\bibitem{liu1} X. Liu and E. S. Titi, Local well-posedness of strong solutions to the three-dimensional compressible Primitive Equations, {\em arxiv1806.09868v1}.

\bibitem{liu2} X. Liu and E. S. Titi, Global existence of weak solutions to the compressible Primitive Equations of atmosphereic dynamics with degenerate viscositites, {\em SIAM J. Math. Anal.}, 51 (2019), 1913-1964.

\bibitem{liu3} X. Liu and E. S. Titi, Zero Mach number limit of the compressible Primitive Equations Part I: well-prepared initial data, {\em arxiv1905.09367v1}.

\bibitem{m} N. Masmoudi and T. K. Wong, On the $H^s$ theory of hydrostatic Euler equations, {\em Arch. Ration. Mech. Anal.}, 204 (2012), 231-271.

\bibitem{s} T. Seng\"{u}l and S. H. Wang, Dynamic transitions and baroclinic instability for 3D continuously stratified Boussinesq flows, {\em J. Math. Fluid Mech.}, 20 (2018), 1173-1193.

\bibitem{tang} T. Tang and H. J. Gao, On the stability of weak solution for compressible primitive equations, {\em Acta Appl. Math.}, 140 (2015), 133-145.


\bibitem{t} R. Temam and M. Ziane, {\em Some mathematical problems in geophysical fluid dynamics}, Handbook of Mathematical Fluid Dynamics, 2004.

\bibitem{v} A. Vasseur and C. Yu, Existence of global weak solutions for 3D degenerate compressible Navier-Stokes equations, {\em Invent. Math.}, 206 (2016), 935-974.


\bibitem{w} F. C. Wang, C. S. Dou and Q. S. Jiu, Global weak solutions to 3D compressible primitive equations with density-dependent viscosity, {\em arxiv:1712.04180v1}.

\bibitem{ws}  S. H. Wang and P. Yang, Remarks on the Rayleigh-Benard convection on spherical shells, {\em J. Math. Fluid Mech.}, 15 (2013), 537-552.

\end{thebibliography}
\end{document}